\documentclass[11pt]{article}
\usepackage{latexsym}
\usepackage{color}
\usepackage{amssymb, graphicx, amsmath}
\usepackage{booktabs}
\usepackage{cases}
\usepackage{subfigure}
\usepackage{graphicx}
\usepackage{CJK,graphicx}

\newcommand{\nn}{\nonumber}

\definecolor{myblue}{rgb}{0,0,0.5}
\definecolor{mygreen}{rgb}{0,0.5,0}
\definecolor{myred}{rgb}{0.5,0,0}
\def \u{{\mbox{\boldmath $u$}}}

\def \w{{\mbox{$\boldsymbol{w}$}}}

\newtheorem{theorem}{Theorem}[section]

\newtheorem{remark}[theorem]{Remark}

\makeatletter
\def\hlinew#1{%
  \noalign{\ifnum0=`}\fi\hrule \@height #1 \futurelet
   \reserved@a\@xhline}

\def \[{\begin{equation}}
\def \]{\end{equation}}
\begin{document}
\begin{CJK*}{GBK}{song}

\begin{center}

{\large \bf  A Revisit of Chen-Teboulle's  Proximal-based Decomposition Method
}\\

\bigskip
\medskip

 {\bf Feng Ma}\footnote{\parbox[t]{16.0cm}{
High-Tech Institute of Xi'an, Xi'an, 710025, Shaanxi, China. This author was supported by the NSFC Grant
 11701564  and
the NSF of Shaanxi Province Grant 2020JQ-485. Email:
 mafengnju@gmail.com}}

\medskip

\today

\end{center}

\bigskip

{\small

\parbox{0.95\hsize}{

\hrule

\medskip

{\bf Abstract.} In this paper, we show that Chen-Teboulle's proximal-based decomposition method can be interpreted as a proximal augmented Lagrangian method. More precisely, it coincides with a  linearized  augmented Lagrangian method. We then proposed three generalized methods based on this interpretation. By invoking recent work (He et
al., IMA J. Numer. Anal., 32 (2020), pp. 227--245), we show that the step size condition of Chen-Teboulle's method  can be relaxed without adding any further assumptions. Our analysis offers a new insight into this proximal-based decomposition method.
\medskip


 \medskip

  \hrule

  }}

\bigskip

\section{Introduction}

In this paper, we consider the following convex problem with separable structure
\[  \label{Problem-LC}
\min \{f(x)+g(z): A x=z\}
  \]
  where $A\in \mathbb{R}^{m \times n}$, $f: \mathbb{R}^{n} \mapsto(-\infty,+\infty)$ and $g: \mathbb{R}^{m} \mapsto(-\infty,+\infty)$ are given closed proper convex functions. The solution set of (\ref{Problem-LC}) is assumed to be nonempty.

The problem   \eqref{Problem-LC} was  discussed in the early days  and has received a lot of attentions in recent years
since it  captures many practical models in sparse and low-rank optimization. Some classical algorithms such as augmented Lagrangian method (ALM) \cite{Hes69,Powell69}, alternating direction method of multipliers (ADMM) \cite{Gab83,GM}, proximal gradient method \cite{Nesterov,Beck2009}   have gained much popularity due to their flexibility and efficiency for solving \eqref{Problem-LC}. In this paper,   we restrict our attention to   Chen-Teboulle's proximal-based decomposition method \cite{Chen94}, which is also  well-known to handle separable convex optimization problems.
The method was dated back from the seminal  work of Chen and Teboulle  in early 1990s  and
soon became a popular distributed primal dual algorithm in convex programming.
Differ from  the ADMM that combines
both ideas of the augmented Lagrangian framework and alternating strategy,  Chen and Teboulle's algorithm  \cite{Chen94} exploits directly for the Lagrangian function \eqref{Aug-L-F} and enjoys parallel architecture. More specially, let the Lagrangian function associated with (\ref{Problem-LC}) be defined by
  \[  \label{Aug-L-F}
L(x, z, y):=f(x)+g(z)+ \left\langle y, A x-z\right\rangle,
    \]
where $y \in \mathbb{R}^{m}$ be the Lagrange multiplier associated  with the linear constraint in (\ref{Problem-LC}).  Chen-Teboulle's algorithm  \cite{Chen94} iterates as follows:
   \begin{subequations} \label{CT-1}
  \begin{numcases}{\hbox{\quad}}
  \label{CT-1-p}  p^{k+1}=\arg \max \left\{L\left(x^{k}, z^{k}, y\right)-\left(1 /\left(2 \lambda\right)\right)\left\|y-y^{k}\right\|^{2}\right\},\\[0.2cm]
 \label{CT-1-x}  x^{k+1}=\arg \min \left\{L\left(x, z^{k}, p^{k+1}\right)+\left(1 /\left(2 \lambda\right)\right)\left\|x-x^{k}\right\|^{2}\right\},\\[0.2cm]
 \label{CT-1-z} z^{k+1}=\arg \min \left\{L\left(x^{k}, z, p^{k+1}\right)+\left(1 /\left(2 \lambda\right)\right)\left\|z-z^{k}\right\|^{2}\right\},\\[0.2cm]
 \label{CT-1-y} y^{k+1}=\arg \max \left\{L\left(x^{k+1}, z^{k+1}, y\right)-\left(1 /\left(2 \lambda\right)\right)\left\|y-y^{k}\right\|^{2}\right\},
\end{numcases}
\end{subequations}
where $\lambda>0$ is  a proximal  parameter. For simplicity, we fix it in our discussion.
From  \eqref{CT-1}, we can observe that the algorithm performs two proximal steps in the dual variables and
one proximal step for the primal variables  $x$ and $z$ separately; and the multiplier variable is first predicted by $p^{k+1}$ step \eqref{CT-1-p}  and then corrected by the  $y^{k+1}$  \eqref{CT-1-y}. With these structural
features, the algorithm \eqref{CT-1} is also called a predictor
corrector proximal multiplier method.

Recall the Lagrangian function defined by \eqref{Aug-L-F}, the  scheme  \eqref{CT-1} can  be rewritten as
   \begin{subequations} \label{CT-2}
  \begin{numcases}{\hskip-1cm\hbox{(PCPM)}}
  \label{CT-2-p}  p^{k+1}=y^{k}+\lambda\left(A x^{k}-z^{k}\right),\\[0.2cm]
 \label{CT-2-x}  x^{k+1}=\arg \min \left\{f(x)+\left\langle p^{k+1}, A x\right\rangle+\left(1 /\left(2 \lambda\right)\right)\left\|x-x^{k}\right\|^{2}\right\},\\[0.2cm]
 \label{CT-2-z} z^{k+1}=\arg \min \left\{g(z)-\left\langle p^{k+1}, z\right\rangle+\left(1 /\left(2 \lambda\right)\right)\left\|z-z^{k}\right\|^{2}\right\},\\[0.2cm]
 \label{CT-2-y} y^{k+1}=y^{k}+\lambda\left(A x^{k+1}-z^{k+1}\right),
\end{numcases}
\end{subequations}
The PCPM only involves evaluations of proximal operators  which are for many problems in closed-form or simple to compute. As a consequence, the PCPM's subproblems may be much easier to solve than  other methods based on the augmented Lagrangian function  such as  ALM, ADMM. Since  the primal subproblems \eqref{CT-1-x}  and \eqref{CT-1-z} are two
separate minimizations, this algorithm is also suitable for parallel computation.  These features make the implementation of PCPM very easy. In \cite{Chen94},   it is proved that PCPM converges to the solutions of \eqref{Problem-LC}
when
\[\label{CT-size1}
\lambda \leq \frac{1}{2\max(\|A\|,1) },
\]
where $\|\cdot\|$ is the spectral norm.

Recently, there are several algorithmic frameworks are developed for unifying and analyzing  a class of primal dual algorithms for the problem  \eqref{Problem-LC}. These frameworks encompass   the PCPM  as special cases.
 For example, it was shown in \cite{Solodov}  that PCPM  can be regarded as a special case in the hybrid inexact proximal point framework developed
in \cite{Solodov1999}. In \cite{He2012}, the PCPM is categorized into a unified framework of proximal-based decomposition methods for monotone variational inequalities. In  \cite{Deng2017}  it is shown that the PCPM is equivalent to a Jacobi-Proximal ADMM.  The convergence results of PCPM can be reconducted by  all these interpretations and the condition \eqref{CT-size1} can also be recovered.

As delineated in the literature \cite{Chen94,Shefi2014,Becker}, the parameter $\lambda$ determines the step size for solving the  subproblems, and it is important to choose
appropriate value for  $\lambda$  to ensure PCPM's efficiency. Intuitively, if $\lambda$  is larger,  each proximal term  in the  subproblems could play a lighter weight in the objective and thus the  variables can be updated with  larger step sizes. If  $\lambda$  is tiny, then it implies that the  subproblems are solved conservatively with a too-small step size;  and it might be not  preferable from numerical perspective. Therefore,  it would be desirable to  consider the possibility of further relaxing the condition \eqref{CT-size1} as long as the convergence of (\ref{LALM}) can be guaranteed, so that larger values of $\lambda$ can be chosen. If this case is possible, one can expect a further speedup of the convergence of PCPM \eqref{CT-2} without additional computation. In  \cite{Shefi2014}, Shefi and  Teboulle showed that
PCPM can be viewed as a linearization of the quadratic penalty term in the parallel decomposition of the proximal method of multipliers and improves the condition to
\[\label{CT-size7}
\lambda \leq \frac{1}{\sqrt{2}\max(\|A\|,1) }.
\]
In the work \cite{Ma14},  by following the same analysis as in \cite{He2012},  this condition is further relaxed to
\[\label{CT-size2}
\lambda < 1 / \sqrt{\|A\|^{2}+1},
\]
which is less restrictive than the conditions \eqref{CT-size1} and \eqref{CT-size7}.
In very recent work \cite{Becker}, Becker demonstrates that the PCPM can be viewed as a preconditioned  proximal point algorithm applied to the primal-dual formulation of problem \eqref{Problem-LC}, which also gets the same condition \eqref{CT-size2}.

Based on above results, it is natural to ask whether the condition \eqref{CT-size2} for   the PCPM  is optimal enough to ensure the the convergence, or more precisely, necessary and sufficient  to ensure the the convergence.  To answer this question,
 we show that, the PCPM can be interpreted as a proximal ALM.  As we will see in Section \ref{sect-2}, the PCPM algorithm
coincides with the linearized ALM algorithm applied to a block reformulation of \eqref{Problem-LC}.   This interpretation has
interesting implications for the analysis of the PCPM  algorithm. First, some known convergence results of PCPM can be easily  recovered or even simplified from the proximal ALM.  Second, it also allows us to translate some useful variants of  the  proximal ALM from the literature  to the PCPM algorithm. We shall provide three generalized PCPM algorithms for the problem \eqref{Problem-LC}.
Then, motivated by recent work of  indefinite proximal ALM \cite{HMY-ALM}, we  show that the convergence condition of  PCPM  can be further relaxed by
\[\label{CT-size3}
 \lambda< \frac{1}{\sqrt{\frac{3}{4}(\|A\|^2+1)}},
\]
without making further assumptions.

The paper is organized as follows. In Section \ref{sect-2},  we briefly review the proximal ALM, and establish the connection
between proximal ALM and PCPM. Then,   in Section \ref{Sec:bound},  we focus on a general PCPM and discuss its optimal condition bound.
In Section \ref{sect-Ext},  we present a more general  extension  of the PCPM with different step sizes for the  linearly constrained convex minimization model. Finally, some conclusions are drawn in Section \ref{Sec:conclusion}.

 \section{  The PCPM is a Proximal ALM}\label{sect-2}

\setcounter{equation}{0}
 In this section, we show that the PCPM \eqref{CT-2} is a special case of the proximal  augmented Lagrangian method  with a particular proximal regularization term. The convergence condition  \eqref{CT-2} can be easily rediscovered by this interpretation.

 \subsection{Proximal ALM}
For convenience of our demonstration, we first introduce some auxiliary variables and   reformulate the problem  into a block version.

Let
\[ \label{VI-P2-wF}
            \w = \left(\begin{array}{c}
                     x\\  z\\ y\end{array} \right), \quad
                          \u = \left(\begin{array}{c}
                     x\\
                 z\end{array} \right),      \quad  M=(A,-I),
 \quad \hbox{and}    \quad {\boldsymbol{\theta}}(\u)=f(x) +
                     g(z).
   \]
   Then the problem \eqref{Problem-LC} can be
rewritten as
\[  \label{Problem-LC2}
\min \{\theta(\u): M\u=0\}.
  \]
The augmented Lagrangian function associated with problem (\ref{Problem-LC2}) is given by
\[  \label{Aug-L-F2}
   {\cal L}_{\lambda}(\u,y) ={\boldsymbol{\theta}}(\u)+\left\langle y, M\u\right\rangle+\left(\lambda /2 \right)\left\|M\u\right\|^{2},
    \]
with $\lambda>0$ the penalty parameter for the linear constraints. Given a starting vector $\left(\u^{0}, y^{0}\right) \in \mathbb{R}^{n} \times \mathbb{R}^{m} \times \mathbb{R}^{m}$, the augmented Lagrangian method (ALM) originally proposed in \cite{Hes69,Powell69} for (\ref{Problem-LC}) generates iterations as
\begin{subequations} \label{ALM}
\begin{numcases}{\hbox{(ALM)\quad}}
\label{ALM-x}  \u^{k+1}=\arg\min \bigl\{ {\cal L}_{\lambda}(\u,y^k)\; \big| \; \u\in\mathbb{R}^{n \times m}
    \bigr\},\\
\label{ALM-l} y^{k+1}=y^{k}+\lambda \left(M\u^{k+1}\right).
\end{numcases}
\end{subequations}
The  computational complexity of the ALM algorithm is dominated by the primal subproblem,  so it is meaningful to discuss how to efficiently solve  (\ref{ALM-x}). An interesting strategy is to regularize the primal subproblem (\ref{ALM-x}) by a quadratic proximal term and accordingly get the proximal version of ALM:
\begin{subequations} \label{PALM}
\begin{numcases}{\hbox{(Proximal ALM)\quad}}
\label{PALM-x}  \u^{k+1} =\arg\min \bigl\{{\cal L}_{\lambda}(\u,y^k) +\frac{1}{2}\|\u-\u^k\|_{\cal P}^2\; \big| \; \u\in\mathbb{R}^{n \times m}\bigr\},\\
\label{PALM-l} y^{k+1}=y^{k}+\lambda \left(M\u^{k+1}\right).
\end{numcases}
\end{subequations}
In (\ref{PALM-x}), $\frac{1}{2}\|\u-\u^k\|_{\cal P}^2$ is the quadratic proximal regularization term and ${\cal P}$ is the proximal matrix that is usually required to be positive definite in the literature.

Typically, we can  linearize the augmented term by choosing appropriate ${\cal P}$  as
 \[\label{ALM-P}
 {\cal P} = \eta I-\lambda M^{\top}M.
 \]
 At this case, the primal subproblem (\ref{PALM-x}) is specified as
\[\label{ALM-x-4}
 \u^{k+1} =\arg\min \bigl\{  {\boldsymbol{\theta}}(\u) +\frac{\eta }{2}\|\u - \u^k-\frac{1}{\eta }M^{\top} \bigl(y^k +\lambda(M\u^k)\big)\|^2\; \big| \; \u\in\mathbb{R}^{n \times m}\},
 \]
which amounts to estimating the proximity operator of ${\boldsymbol{\theta}}(\u)$. The implementation for such cases is usually simple.

Hence, the linearized ALM, which is a special case of the proximal ALM (\ref{PALM}) with ${\cal P}$  given in \eqref{ALM-P}, reads as
\begin{subequations} \label{LALM}
\begin{numcases}{\hskip-1cm\hbox{(Linearized ALM)}}
\label{LALM-x}   \u^{k+1} =\arg\min \bigl\{  {\boldsymbol{\theta}}(\u) +\frac{\eta }{2}\|\u - \u^k-\frac{1}{\eta }M^{\top} \bigl(y^k +\lambda(M\u^k)\big)\|^2\; \big| \; \u\in\mathbb{R}^{n \times m }\},\ \\
   \label{LALM-l} y^{k+1}=y^{k}+\lambda \left(M\u^{k+1}\right).
\end{numcases}
\end{subequations}
For the linearized ALM (\ref{LALM}) in the literature, the parameter $\eta$ is required to satisfy the condition $\eta>\lambda \|M\|$ so as to ensure the positive definiteness of the matrix ${\cal P}$ given in (\ref{ALM-P}) and hence the convergence of (\ref{LALM}). We refer to \cite{YangYuan,HMY-ALM} for the detail of convergence analysis of the linearized ALM (\ref{LALM}).

 \subsection{Proximal  ALM Perspective}
In this subsection, we show the  PCPM \eqref{CT-2} is a special case of the proximal ALM \eqref{PALM}. This will be done by simple algebraic manipulation and simplification.

First, substituting    $p^{k+1}$ \eqref{CT-2-p}  into the update for $x^{k+1}$ \eqref{CT-2-x} and $z^{k+1}$ \eqref{CT-2-z},
the  primal iterations are
\[
  \label{CT-3-p}
 x^{k+1}=\arg \min \left\{f(x)+ \left\langle y^{k}+\lambda\left(A x^{k}-z^{k}\right), Ax\right\rangle+\left(1 /\left(2 \lambda\right)\right)\left\|x-x^{k}\right\|^{2}\right\},
 \]
 and
 \[
 \label{CT-3-z} z^{k+1}=\arg \min \left\{g(z)- \left\langle y^{k}+\lambda\left(A x^{k}-z^{k}\right), z\right\rangle +\left(1 /\left(2 \lambda\right)\right)\left\|z-z^{k}\right\|^{2}\right\},
\]
respectively. Note that the update $p$ is eliminated. We  then have
 \begin{eqnarray} \label{CT-3-x-z1}
( x^{k+1},z^{k+1})=&\arg \min &\left\{f(x)+g(z)+\left\langle y^{k}+\lambda\left(A x^{k}-z^{k}\right), Ax-z\right\rangle \right.  \nn\\
&& \left.+\left(1 /\left(2 \lambda\right)\right)\left\|x-x^{k}\right\|^{2}+\left(1 /\left(2 \lambda\right)\right)\left\|z-z^{k}\right\|^{2}\right\}\nn\\
=&\arg \min &\left\{f(x)+g(z)+\left\langle y^{k}, Ax-z\right\rangle+\left(\lambda /2 \right)\left\|Ax-z\right\|^{2} \right.  \nn\\
&& +\left(1 /\left(2 \lambda\right)\right)\left\|x-x^{k}\right\|^{2}+\left(1 /\left(2 \lambda\right)\right)\left\|z-z^{k}\right\|^{2} \nn\\
&& \left.-\left(\lambda /2 \right)\left\|Ax-z\right\|^{2}+\left\langle \lambda \left(A x^{k}-z^{k}\right),  Ax-z \right\rangle  \right\}. \nn\\
=&\arg \min &\left\{f(x)+g(z)+\left\langle y^{k}, Ax-z\right\rangle+\left(\lambda /2 \right)\left\|Ax-z\right\|^{2} \right.  \nn\\
&& +\left(1 /\left(2 \lambda\right)\right)\left\|x-x^{k}\right\|^{2}+\left(1 /\left(2 \lambda\right)\right)\left\|z-z^{k}\right\|^{2} \nn\\
&& \left.-\left(\lambda /2 \right)\left\| \left(Ax-z\right)- \left(A x^{k}-z^{k}\right)\right\|^{2}+\left(\lambda /2 \right)\left\|Ax^{k}-z^{k}\right\|^{2}  \right\}.
   \end{eqnarray}
  Ignoring
some constant terms in the minimization problem of the last equality, we have
 \begin{eqnarray} \label{CT-3-x-z2}
( x^{k+1},z^{k+1})=&\arg \min &\left\{f(x)+g(z)+\left\langle y^{k}, Ax-z\right\rangle+\left(\lambda /2 \right)\left\|Ax-z\right\|^{2} \right.  \nn\\
&& +\left(1 /\left(2 \lambda\right)\right)\left\|x-x^{k}\right\|^{2}+\left(1 /\left(2 \lambda\right)\right)\left\|z-z^{k}\right\|^{2} \nn\\
&& \left.-\left(\lambda /2 \right)\left\| A\left(x-x^{k}\right)- \left(z-z^{k}\right)\right\|^{2}  \right\}.
   \end{eqnarray}
The above scheme can be represented as
         \begin{subequations} \label{CT-3-x-z}
    \[ \label{CT-3-x-z3}
( x^{k+1},z^{k+1})=\arg \min \left\{f(x)+g(z)+\left\langle y^{k}, Ax-z\right\rangle+\left(\lambda /2 \right)\left\|Ax-z\right\|^{2} +\frac{1}{2}\left\|\left(\begin{array}{c} x-x^{k}\\ z-z^{k}\end{array}\right)\right\|_{\cal P}^{2}  \right\},\]
with
\[\label{CT-P}
{\cal P}=\left(\begin{array}{cc}
     \frac{1}{\lambda}I- \lambda A^{\top}A     &    \lambda A^{\top} \\
               \lambda A &     (\frac{1}{\lambda}-\lambda)I  \end{array} \right).
\]
   \end{subequations}
Recall  \eqref{PALM} and the definitions in  \eqref{VI-P2-wF}, we can see the algorithm consits of \eqref{CT-3-x-z} and \eqref{CT-2-y}
 is a special proximal ALM.

To summarize,  the PCPM \eqref{CT-2}   is interpreted as a  proximal ALM applied to problem \eqref{Problem-LC2}. This is different from the result in \cite{Becker},
in which the algorithm is interpreted as the preconditioned proximal point algorithm applied to a primal-dual
reformulation of the original problem \eqref{Problem-LC}.

\begin{remark}
Let us take a deeper look at the regularization matrix ${\cal P}$  in \eqref{CT-P} and derive a simpler representation of it. We have
\begin{eqnarray} \label{Matrix-G-P}
{\cal P} &=   &\left(\begin{array}{cc}
     \frac{1}{\lambda}I_n- \lambda A^{\top}A     &    \lambda A^{\top} \\
               \lambda A &     (\frac{1}{\lambda}-\lambda)I_m  \end{array} \right) \nn \\
         &=   &  \frac{1}{\lambda}I-\lambda\cdot \left(\begin{array}{cc}
   A^{\top}A     &      -A^{\top} \\
             - A &     I  \end{array} \right) \nn \\
                      & \stackrel{\eqref{VI-P2-wF}}{=}   &  \frac{1}{\lambda}I-\lambda \cdot M^{\top}M.
    \end{eqnarray}
\end{remark}
We can see that here ${\cal P}$ is a special case of \eqref{ALM-P} where $\eta=  \frac{1}{\lambda}.$ Hence, the PCPM \eqref{CT-2}  can be further  interpreted as a  linearized ALM with ${\cal P}$  given by \eqref{CT-P}.

\begin{remark}
The step size condition \eqref{CT-size2} can be easily obtained by this proximal ALM's perspective. Since the proximal regularization matrix ${\cal P}$  is usually required to be positive definite, we have
\[\label{CT-P2}
{\cal P}=\frac{1}{\lambda}I-\lambda \cdot M^{\top}M \succ 0.
\]
It only remains to ensure
\[
 \lambda^2 <\frac{1}{\|M^{\top}M\|}.
\]
Recall $M$ defined in \eqref{VI-P2-wF}, the condition reduces to $\lambda^2(\|A\|^2+1)<1$. Thus the step size condition \eqref{CT-size2} is obtained.
\end{remark}

 \subsection{Two Generalized PCPM}
 In this section, we present two generalized PCPM. The first one is developed by following the idea of relaxing ALM. The second  is developed by viewing it as a variant of PPA.

 In order to accelerate the convergence of ALM or the proximal ALM, one practical strategy is to attach  a relaxation
factor to the Lagrange-multiplier-updating step in the algorithm. For the proximal ALM, the relaxed scheme is
 \begin{subequations} \label{PALM2}
\begin{numcases}{\hbox{\quad}}
\label{PALM2-x}  \u^{k+1} =\arg\min \bigl\{{\cal L}_{\lambda}(\u,y^k) +\frac{1}{2}\|\u-\u^k\|_{\cal P}^2\; \big| \; \u\in\mathbb{R}^{n \times m}
    \bigr\},\\
\label{PALM2-l} y^{k+1}=y^{k}+\gamma\lambda \left(M\u^{k+1}\right).
\end{numcases}
\end{subequations}
 where the relaxation
factor $\gamma$ can be chosen in the interval $(0,2)$, Recall that the proximal ALM (\ref{PALM}) is a special case of (\ref{PALM2}) with $\gamma=1$. Numerically, an overrelaxation choice $\gamma\in [1.5, 1.8]$ can usually lead to faster convergence; see some
numerical results in \cite{Ma2019}.
Since Chen-Teboulle's algorithm is a proximal ALM, we can relax its dual step size as the proximal ALM  and get the following relaxed algorithm.
   \begin{subequations} \label{CT-4}
  \begin{numcases}{\hskip-1cm\hbox{(G-PCPM-I)}}
  \label{CT-4-p}  p^{k+1}=y^{k}+\lambda\left(A x^{k}-z^{k}\right),\\[0.2cm]
 \label{CT-4-x}  x^{k+1}=\arg \min \left\{f(x)+\left\langle p^{k+1}, A x\right\rangle+\left(1 /\left(2 \lambda\right)\right)\left\|x-x^{k}\right\|^{2}\right\},\\[0.2cm]
 \label{CT-4-z} z^{k+1}=\arg \min \left\{g(z)-\left\langle p^{k+1}, z\right\rangle+\left(1 /\left(2 \lambda\right)\right)\left\|z-z^{k}\right\|^{2}\right\},\\[0.2cm]
 \label{CT-4-y} y^{k+1}=y^{k}+\gamma\lambda\left(A x^{k+1}-z^{k+1}\right),
\end{numcases}
\end{subequations}
where $\gamma\in(0,2)$.

In the PPA literature, it is commonly known that  the PPA scheme can be  relaxed, i.e., we can generate the new iterate by relaxing the output of the original PPA  appropriately.  This is usually based on combining the output of the operation with the former iterate. On the other hand, the proximal ALM can be interpreted as a type of preconditioned proximal point algorithm (PPA), see \cite{Gu2014} for details. We  refer to   \cite{Becker}  for direct
discussions on the interpretation.  Hence, the PCPM, as a special PPA, can also be generalized. More
precisely, let the output point of \eqref{CT-2} be denoted by $\tilde{\w}^k$, then the relaxed PCPM yields the
new iterate via
\begin{subequations}  \label{PSALM-A}
\begin{numcases}{\hskip-1cm\hbox{(G-PCPM-II)}}
   \nonumber\\[-0.1cm]
 \label{PSALMA-X}
    \! \left\{\begin{array}{l}
  \label{CT-5-p}  p^{k+1}=y^{k}+\lambda\left(A x^{k}-z^{k}\right),\\[0.2cm]
 \label{CT-5-x}  \tilde{x}^{k}=\arg \min \left\{f(x)+\left\langle p^{k+1}, A x\right\rangle+\left(1 /\left(2 \lambda\right)\right)\left\|x-x^{k}\right\|^{2}\right\},\\[0.2cm]
 \label{CT-5-z}  \tilde{z}^{k}=\arg \min \left\{g(z)-\left\langle p^{k+1}, z\right\rangle+\left(1 /\left(2 \lambda\right)\right)\left\|z-z^{k}\right\|^{2}\right\},\\[0.2cm]
 \label{CT-5-y}  \tilde{y}^{k}=y^{k}+\lambda\left(A x^{k+1}-z^{k+1}\right),
    \end{array} \right. \\[0.2cm]
\label{PSALMA-L}
    \;\;   {\w}^{k+1} = {\w}^k - \gamma\Bigl({\w}^k- \tilde{\w}^{k}\Bigr),
    \end{numcases}
\end{subequations}
where $\gamma \in(0,2)$ is the relaxation factor. In particular, $\gamma$ is called an under-relaxation  factor when $\gamma \in(0,1)$  or over-relaxation factor when   $\gamma \in(1,2)$; and the relaxed G-PCPM-II \eqref{PSALM-A} reduces to the original PCPM \eqref{CT-2} when $\gamma=1$.
\section{Optimal bound  on Step Sizes}\label{Sec:bound}

\setcounter{equation}{0}
In this section, with the proximal ALM interpretation, we shall show that
 the step size parameters $\lambda$ and $\gamma$ in the generalized PCPM \eqref{CT-4} can be related by the formula
\[\label{new-bound}
 \lambda<  \frac{1}{\sqrt{\frac{2+\gamma}{4}(\|A\|^2+1)}}.
\]
Note that when $\gamma=1$, G-PCPM-I \eqref{CT-4} reduces to PCPM and the above condition reduces to $\lambda< \frac{1}{\sqrt{\frac{3}{4}(\|A\|^2+1)}}$ which improves the results  in the works \cite{Chen94,Ma14,Becker} to ensure the convergence.

The new bound \eqref{new-bound} relies on  the convergence results studied in \cite{HMY-ALM}  for the proximal ALM. Here, we describe
the main results for the proximal ALM in \cite{HMY-ALM} by the following presentation,
but omit the proof.
\begin{center}
\fbox{
\begin{minipage}{16cm}
In \cite{HMY-ALM}, the authors showed that for the proximal ALM
\begin{subequations} \label{LP-ALM}
\begin{numcases}{\hbox{(IDP-ALM)\quad}}
\label{LP-ALM-x}  \u^{k+1} =\arg\min \bigl\{{\cal L}_{\lambda}(\u,\lambda^k) +
    \frac{1}{2}\|\u-\u^k\|_{{\cal P}}^2 \; \big| \; \u\in\mathbb{R}^{l}
    \bigr\},\\
\label{LP-ALM-l}y^{k+1} = y^k +\gamma \lambda(M\u^{k+1}), \quad \gamma \in (0,2),
\end{numcases}
\end{subequations}
the proximal matrix ${\cal P}$ in the term $\frac{1}{2}\|\u-\u^k\|_{{\cal P}}^2$  can be indefinite without any further assumptions. In particular, Let ${\cal P}$ be specified by the structure
\begin{subequations}\label{D0-L}
\[
{\cal P}={\cal D}-(1-\tau)\lambda M^{\top}M,
\]
where ${\cal D}$ is an arbitrarily positive definite matrix in $\mathbb{R}^{l}$.  Then when
\[\label{taugamma}
 \tau > \frac{2+\gamma}{4},
 \]
 \end{subequations}
the IDP-ALM \eqref{LP-ALM} converges globally to a
solution of \eqref{Problem-LC2}.
\end{minipage} }
\end{center}

Since the proximal ALM contains PCPM as a special case, the proximal matrix ${\cal P}$ defined in \eqref{CT-3-x-z} can  employ indefinite setting according to above results. To do this, we just need to choose $\lambda $ to guarantee
  \[\label{D0-L2}
{\cal D}={\cal P}+(1-\tau)\lambda M^{\top}M\succ 0, \quad \hbox{for}\; \tau\in \Bigl(\displaystyle\frac{2+\gamma}{4}, 1\Bigr).
\]
To fulfill (\ref{D0-L2}), notice that
\begin{eqnarray}  \label{D0-L1}
{\cal D}&\stackrel{\eqref{Matrix-G-P}}{=} &\frac{1}{\lambda}I-\lambda \cdot M^{\top}M+ (1-\tau)\lambda M^{\top}M\\
&=&\frac{1}{\lambda}I-\tau\lambda M^{\top}M.
  \end{eqnarray}
Recall $M$ defined in \eqref{VI-P2-wF}. Then we have
\[\label{CT-new-size1}
\lambda<  \frac{1}{\sqrt{\tau(\|A\|^2+1)}}  \quad \Rightarrow     \quad   {\cal D} \succ 0.
\]
Note that $\tau\in \Bigl(\displaystyle\frac{2+\gamma}{4}, 1\Bigr)$ is arbitrary, we have
\[\label{CT-new-size2}
 \lambda< \frac{1}{\sqrt{\frac{2+\gamma}{4}(\|A\|^2+1)}}.
\]
In Figure 1, we plot the evolutions of the step size $\lambda$  with respect
to the norm $\|A\|$ for the three stepsize conditions of the PCPM.  The ratios of  the two step sizes \eqref{CT-size2}   and  \eqref{CT-size3} to  \eqref{CT-size1}  are displayed in Figure 2. These plots show that the stepsize condition is enlarged, and a larger value of  $\lambda$ seems more preferable in practice because it can yield
a larger step size .

\begin{figure}
\centering
\begin{minipage}[c]{.7\textwidth}
\centering
\includegraphics[width =\textwidth]{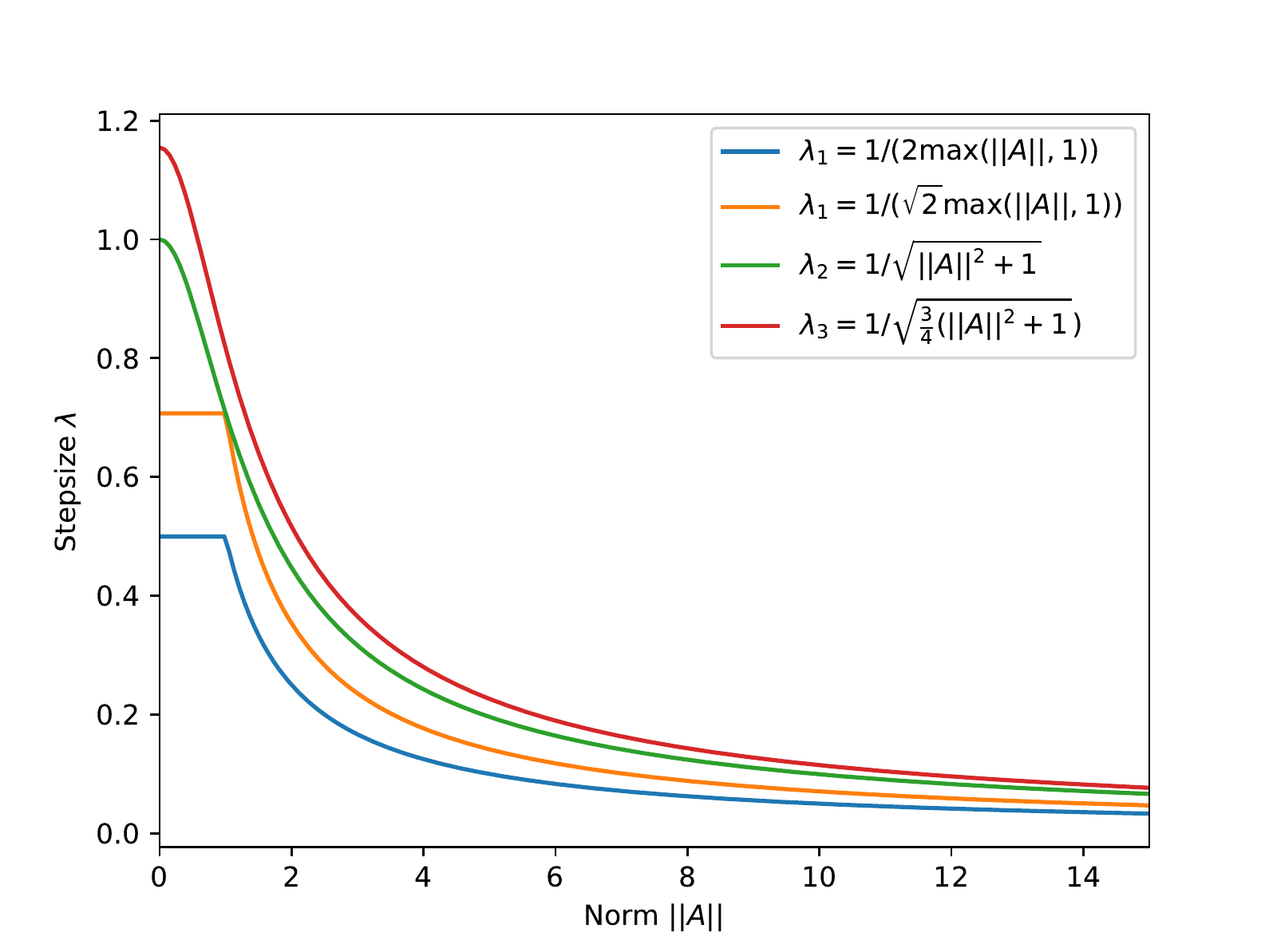} \\
\end{minipage}
\caption{The curves of three different step sizes.} \label{Figure-MCP}
\end{figure}
\begin{figure}
\centering
\begin{minipage}[c]{.7\textwidth}
\centering
\includegraphics[width =\textwidth]{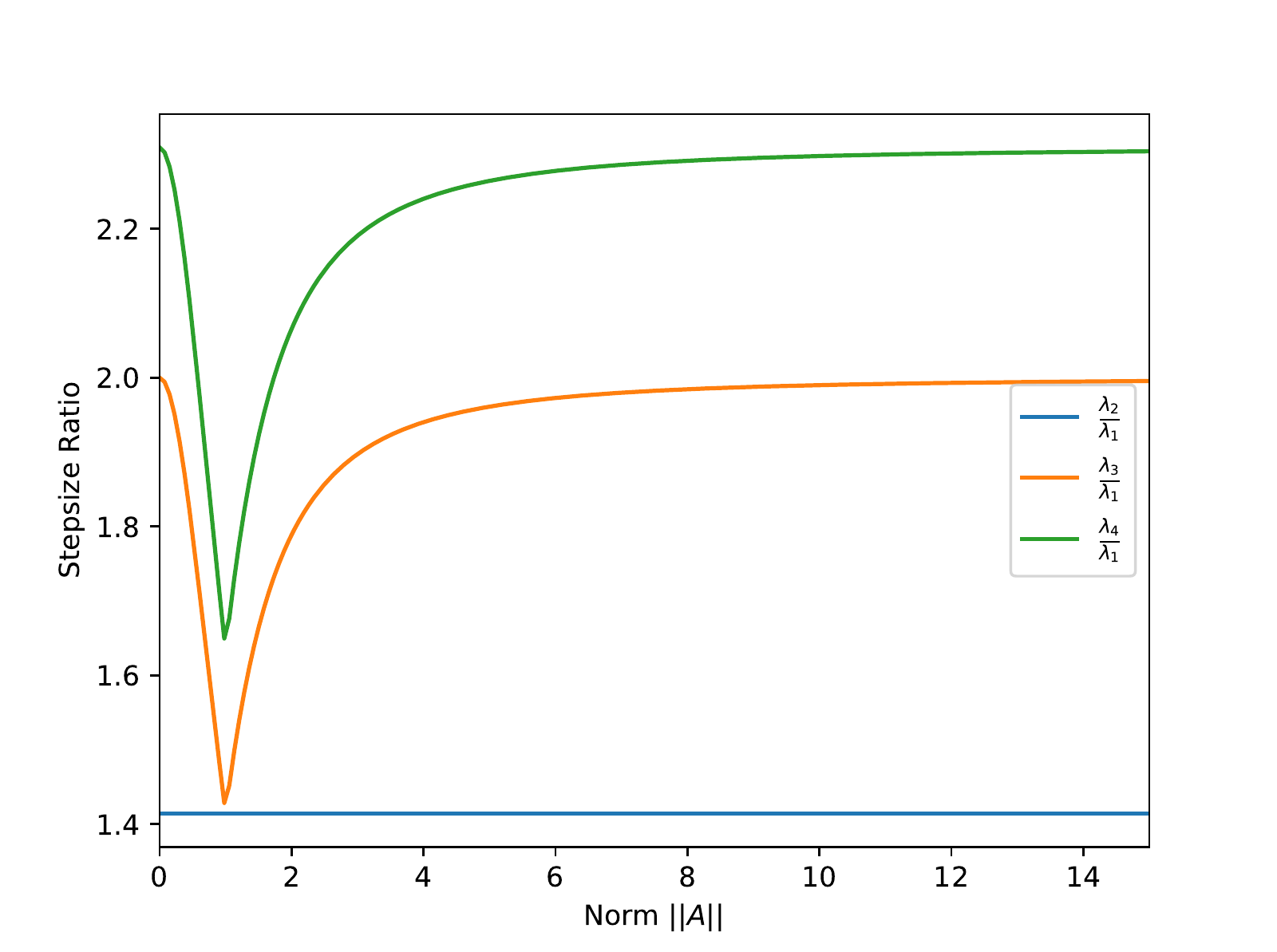} \\
\end{minipage}
\caption{Illustration of the step size ratio.} \label{Figure-MCP2}
\end{figure}

\begin{remark}
In \cite{HMY-ALM}, a counterexample was given showing that the convergence condition \eqref{D0-L} is  optimal for the proximal ALM. Since PCPM is a special case of the proximal ALM, the new bound \eqref{CT-new-size2} thus can not be further improved if we don't add any further assumptions on the model.
\end{remark}

\section{Further Extension}\label{sect-Ext}

\setcounter{equation}{0}
Technically, we can  extend  the original PCPM scheme  and its variants to handle
the following general problem
\[  \label{Problem-LC3}
\min \{f(x)+g(z): A x +Bz=b\},
  \]
    where $A\in \mathbb{R}^{m \times n}, B\in \mathbb{R}^{m \times l}$, $f: \mathbb{R}^{n} \mapsto(-\infty,+\infty)$ and $g: \mathbb{R}^{l} \mapsto(-\infty,+\infty)$ are given closed proper convex functions.

To solve \eqref{Problem-LC3},  instead of  \eqref{CT-2} we propose  the following generalized PCPM
   \begin{subequations} \label{CT-G3}
  \begin{numcases}{\hskip-1cm\hbox{(G-PCPM-III) }}
 p^{k+1}=y^{k}+\lambda\left(A x^{k}+Bz^{k}-b\right),\\[0.2cm]
  x^{k+1}=\arg \min \left\{f(x)+\left\langle p^{k+1}, A x\right\rangle+\left(1 /\left(2 \tau\right)\right)\left\|x-x^{k}\right\|^{2}\right\},\\[0.2cm]
 z^{k+1}=\arg \min \left\{g(z)+\left\langle p^{k+1}, Bz\right\rangle+\left(1 /\left(2 \sigma\right)\right)\left\|z-z^{k}\right\|^{2}\right\},\\[0.2cm]
y^{k+1}=y^{k}+\gamma\lambda\left(A x^{k+1}+Bz^{k+1}-b\right),
\end{numcases}
\end{subequations}
where $\lambda$ is the proximal parameter for the dual regularization; $\tau,\sigma$ are two different positive parameters for the primal regularization; the relaxation factor $\gamma \in(0,2)$.

Now, we discuss how to  drive the step size condition to ensure the convergence of PCPM-III \eqref{CT-G3}. We follow the line of analysis in Section \ref{sect-2}. Let us  transform the G-PCPM-III \eqref{CT-G3}  into an equivalent  proximal ALM with a special matrix ${\cal P}$.
 Ignoring
some constant terms in the minimization problem of the last equality, we have
 \begin{eqnarray} \label{CT-3-x-z5}
( x^{k+1},z^{k+1})=&\arg \min &\left\{f(x)+g(z)+\left\langle y^{k}, A x +Bz-b\right\rangle+\left(\lambda /2 \right)\left\|A x +Bz-b\right\|^{2} \right.  \nn\\
&& +\left(1 /\left(2 \tau\right)\right)\left\|x-x^{k}\right\|^{2}+\left(1 /\left(2 \sigma\right)\right)\left\|z-z^{k}\right\|^{2} \nn\\
&& \left.-\left(\lambda /2 \right)\left\| A\left(x-x^{k}\right)+B\left(z-z^{k}\right)\right\|^{2}  \right\}.
   \end{eqnarray}
The above scheme can be represented as
         \begin{subequations} \label{CT-3-x-6}
\begin{eqnarray}\label{CT-3-x-z7}
( x^{k+1},z^{k+1})=&\arg \min &  \left\{f(x)+g(z)+\left\langle y^{k}, A x +Bz-b\right\rangle\right.
 \nn\\
&&\left.+\left(\lambda /2 \right)\left\|A x +Bz-b\right\|^{2} +\frac{1}{2}\left\|\left(\begin{array}{c} x-x^{k}\\ z-z^{k}\end{array}\right)\right\|_{\cal P}^{2}  \right\},
 \end{eqnarray}
with
\[\label{CT-III-P}
{\cal P}=\left(\begin{array}{cc}
     \frac{1}{\tau}I- \lambda A^{\top}A     &    -\lambda A^{\top}B \\
               -\lambda B^{\top}A &      \frac{1}{\sigma}I-\lambda B^{\top}B  \end{array} \right).
\]
   \end{subequations}
 Hence, the algorithm G-PCPM-III  can be equivalently presented as
   \begin{subequations}
  \begin{numcases}{\hbox{\quad}}
( x^{k+1},z^{k+1})= \arg \min   \left\{f(x)+g(z)+\left\langle y^{k}, A x +Bz-b\right\rangle\right.
 \nn\\
 \qquad\qquad\qquad\qquad\qquad\left.+\left(\lambda /2 \right)\left\|A x +Bz-b\right\|^{2} +\frac{1}{2}\left\|\left(\begin{array}{c} x-x^{k}\\ z-z^{k}\end{array}\right)\right\|_{\cal P}^{2}  \right\}, \\[0.2cm]
 y^{k+1}=y^{k}+\gamma\lambda\left(A x^{k+1}+Bz^{k+1}-b\right),
\end{numcases}
\end{subequations}
 where ${\cal P}$ is given by \eqref{CT-III-P}.

   We can see that
   \begin{eqnarray} \label{Matrix-G-P2}
{\cal P}
         &=   &    \lambda \left(\begin{array}{cc}
     \frac{1}{\tau\lambda }I- A^{\top}A     &    -A^{\top}B \\
                -B^{\top}A &      \frac{1}{\sigma\lambda }I- B^{\top}B  \end{array} \right)\nn \\
      &=   &
       \lambda  \left(\begin{array}{cc}
   \frac{1}{\tau\lambda}I     &      0 \\
             0 &    \frac{1}{\sigma\lambda}I  \end{array} \right)   -\lambda\cdot \left(\begin{array}{cc}
   A^{\top}A     &     A^{\top}B \\
           B^{\top}A &      B^{\top}B \end{array} \right).
    \end{eqnarray}
    Since the proximal regularization matrix ${\cal P}$  is usually required to be positive definite, we need to ensure
\[   \left(\begin{array}{cc}    \frac{1}{\tau\lambda}I     &      0 \\
             0 &    \frac{1}{\sigma\lambda}I  \end{array} \right) \succ    \left(\begin{array}{cc}
   A^{\top}A     &      A^{\top}B \\
            B^{\top}A &      B^{\top}B \end{array} \right)
          \]
   or
   \[   I \succ  \left(\begin{array}{cc}     \sqrt{\tau\lambda}I     &      0 \\
             0 &    \sqrt{\sigma\lambda}I  \end{array} \right)   \left(\begin{array}{cc}
   A^{\top}A     &     A^{\top}B \\
            B^{\top}A &      B^{\top}B \end{array} \right)\left(\begin{array}{cc}     \sqrt{\tau\lambda}I     &      0 \\
             0 &    \sqrt{\sigma\lambda}I  \end{array} \right).
          \]

          Note that $ \left(\begin{array}{cc}
   A^{\top}A     &     A^{\top}B \\
           B^{\top}A &      B^{\top}B \end{array} \right)=\left(\begin{array}{c}
                     A^{\top}\\
                 B^{\top}\end{array} \right)(A, B) $. We just need   to guarantee
  \[    I \succ  \left(\begin{array}{c}
                  \sqrt{\tau\lambda}   A^{\top}\\
                  \sqrt{\sigma\lambda}B^{\top}\end{array} \right)(\sqrt{\tau\lambda} A, \sqrt{\sigma\lambda}B)  .
          \]
So if the step size parameters satisfy
   \[
 \lambda\tau\left\|A^{\top} A\right\| + \lambda\sigma \left\|B^{\top} B\right\| <1,
\]
then the positive definiteness of matrix $P$ is ensured.

 According to the improved convergence result of the proximal ALM \eqref{D0-L}, the proximal matrix ${\cal P}$  can  employ indefinite setting. At this case, the condition is relaxed by
   \[\label{general-condition}
 \lambda\tau\left\|A^{\top} A\right\| + \lambda\sigma \left\|B^{\top} B\right\| <\frac{4}{2+\gamma}.
\]
If all the parameters  $\lambda,\tau,\sigma$  are chosen to be equal, i.e.,
$\lambda=\tau=\sigma$. The resulting condition \eqref{general-condition} reduces to \eqref{CT-new-size2}.
Since asymptotically
this extension has no difference from the PCPM, we skip the  detailed analysis for this scheme.

\section{Concluding Remarks}\label{Sec:conclusion}
In this paper, we study the predictor
corrector proximal multiplier method (PCPM) for convex programming problems, and show  that it is equivalent to a  linearized  augmented Lagrangian method (ALM) with a special regularization term. This interpretation makes it possible to simplify the convergence analysis, and we can  further relax the step size condition of  PCPM by invoking recent improved convergence study of the proximal ALM. It must be mentioned that our result does not rely on any further assumptions of the problems or algorithms.
Since the linearized ALM is extremely popular in recent years, various
variants and theoretical results   have been developed and studied in the literature. Based on our interpretation, these modifications and theoretical results can be easily injected into the PCPM. Thus,
Our analysis builds on the techniques and recent results of   the proximal ALM and gives some insight of PCPM.

\end{CJK*}
\end{document}